\documentclass[journal]{IEEEtran}

\usepackage[latin1]{inputenc}
\usepackage{cite}
\usepackage[dvips]{graphicx}
\usepackage[cmex10]{amsmath}
\usepackage{amsfonts}
\usepackage{algorithmic}
\usepackage{array}
\usepackage[caption=false,font=footnotesize]{subfig}
\usepackage{dblfloatfix}

\usepackage[dvips,colorlinks=false]{hyperref}
\hypersetup{
pdfauthor = {Lars Eirik Danielsen},
pdftitle = {On the classification of Hermitian self-dual additive codes over GF(9)},
pdfkeywords = {self-dual codes, additive codes, codes over GF(9), graph theory, classification, nonbinary quantum codes}
}

\DeclareMathOperator{\GF}{GF}
\DeclareMathOperator{\wt}{wt}
\DeclareMathOperator{\Aut}{Aut}
\DeclareMathOperator{\tr}{Tr}
\DeclareMathOperator{\symp}{Sp}
\DeclareMathOperator{\sym}{Sym}
\DeclareMathOperator{\rank}{rank}

\begin{document}

\title{On the Classification of Hermitian Self-Dual \\ Additive Codes over GF(9)}
\author{Lars Eirik Danielsen%
\thanks{Manuscript created June 13, 2011; revised November 29, 2011.
This research was supported by the Research Council of Norway.
The material in this paper was presented in part at the 2010 IEEE International 
Symposium on Information Theory, Austin, TX, USA, June 2010.
}
\thanks{L.E. Danielsen is with the Department of Informatics,
University of Bergen, N-5020 Bergen, Norway
(e-mail: larsed@ii.uib.no).}
}

\maketitle

\begin{abstract}\boldmath
  Additive codes over $\GF(9)$ that are self-dual with respect to the
  Hermitian trace inner product have a natural application in quantum
  information theory, where they correspond to ternary quantum
  error-correcting codes.  However, these codes have so far received
  far less interest from coding theorists than self-dual additive
  codes over $\GF(4)$, which correspond to binary quantum codes.
  Self-dual additive codes over $\GF(9)$ have been classified up to
  length 8, and in this paper we extend the complete classification to
  codes of length 9 and 10.  The classification is obtained by using a
  new algorithm that combines two graph representations of self-dual
  additive codes.  The search space is first reduced by the fact that
  every code can be mapped to a weighted graph, and a different graph
  is then introduced that transforms the problem of code equivalence
  into a problem of graph isomorphism. By an extension technique, we
  are able to classify all optimal codes of length 11 and 12. There
  are 56\,005\,876 $(11,3^{11},5)$ codes and 6493 $(12,3^{12},6)$
  codes.  We also find the smallest codes with trivial automorphism
  group.
\end{abstract}


\section{Introduction}

\IEEEPARstart{A}{dditive} codes over $\GF(9)$ of length~$n$ are
$\GF(3)$-linear subgroups of $\GF(9)^n$. Such an additive code
contains $3^k$ codewords for some $0 \le k \le 2n$, and is called an
$(n,3^k)$ code. A code $\mathcal{C}$ can be defined by a $k \times n$
\emph{generator matrix} with entries from $\GF(9)$ whose rows span
$\mathcal{C}$ additively.  We denote $\GF(9) = \{0,1,\omega, \omega^2,
\ldots,\omega^7\}$, where $\omega^2 = \omega + 1$.  \emph{Conjugation}
of $x \in \GF(9)$ is defined by $\overline{x} = x^3$.  The \emph{trace
  map}, $\tr : \GF(9) \to \GF(3)$, is defined by $\tr(x) = x +
\overline{x}$.  Following Nebe, Rains, and Sloane~\cite{sdbook}, we
define the \emph{Hermitian trace inner product} of two vectors
$\boldsymbol{u}, \boldsymbol{v} \in \GF(9)^n$ by
\[
(\boldsymbol{u}, \boldsymbol{v}) = \omega^2(\boldsymbol{u} \cdot
\overline{\boldsymbol{v}} - \overline{\boldsymbol{u}} \cdot
\boldsymbol{v}) = \tr(\omega^2 \boldsymbol{u} \cdot
\overline{\boldsymbol{v}}),
\]
where multiplication by $\omega^2$ is necessary because
the skew-symmetric bilinear form $(\boldsymbol{u} \cdot
\overline{\boldsymbol{v}} - \overline{\boldsymbol{u}} \cdot
\boldsymbol{v})$ does not take values in $\GF(3)$~\cite{sdbook}.
We define the \emph{dual} of the
code $\mathcal{C}$ with respect to the Hermitian trace inner product,
$\mathcal{C}^\perp = \{ \boldsymbol{u} \in \GF(9)^n \mid
(\boldsymbol{u}, \boldsymbol{c})=0 \text{ for all } \boldsymbol{c} \in
\mathcal{C} \}$.  $\mathcal{C}$~is \emph{self-orthogonal} if
$\mathcal{C} \subseteq \mathcal{C}^\perp$.  If $\mathcal{C} =
\mathcal{C}^\perp$, then $\mathcal{C}$ is \emph{self-dual} and must be
an $(n,3^n)$ code. 
The class of trace-Hermitian self-dual additive
codes over $\GF(9)$ is also known as $9^{H+}$~\cite{sdbook}.  
The \emph{Hamming weight} of $\boldsymbol{u}$, denoted
$\wt(\boldsymbol{u})$, is the number of non-zero components of
$\boldsymbol{u}$.  The \emph{Hamming distance} between
$\boldsymbol{u}$ and $\boldsymbol{v}$ is $\wt(\boldsymbol{u} -
\boldsymbol{v})$.  The \emph{minimum distance} of the code
$\mathcal{C}$ is the minimal Hamming distance between any two distinct
codewords of $\mathcal{C}$. Since $\mathcal{C}$ is an additive code,
the minimum distance is also given by the smallest non-zero weight of
any codeword in $\mathcal{C}$.  A code with minimum distance~$d$ is
called an $(n,3^k,d)$ code. The \emph{weight distribution} of the code
$\mathcal{C}$ is the sequence $(A_0, A_1, \ldots, A_n)$, where $A_i$
is the number of codewords of weight~$i$.  The \emph{weight
  enumerator} of $\mathcal{C}$ is the polynomial
\[
W(x,y) = \sum_{i=0}^n A_i x^{n-i} y^i,
\]
where we will denote $W(y) = W(1,y)$.  It follows from the
\emph{Singleton bound}~\cite{nonbin1} that any self-dual additive code
must satisfy $d \le \lfloor\frac{n}{2}\rfloor + 1$.  A code is called
\emph{extremal} if it has minimum distance $\lfloor\frac{n}{2}\rfloor + 1$, 
and \emph{near-extremal} if it has minimum distance $\lfloor\frac{n}{2}\rfloor$.  If
a code has the highest possible minimum distance for the given length, it
is called \emph{optimal}. A tighter bound exists for codes over
$\GF(4)$~\cite{calderbank}, but in general the Singleton bound
is the best known upper bound.  Codes that satisfy the Singleton bound
with equality are also known as \emph{maximum distance separable (MDS)
  codes}.  The well-known \emph{MDS conjecture} implies that self-dual
additive MDS codes over $\GF(9)$ must have length $n \le 10$. We have
shown in previous work~\cite{nonbinary} that there are only three
non-trivial MDS codes, with parameters $(4,3^4,3)$, $(6,3^6,4)$, and
$(10,3^{10},6)$, given that the MDS conjecture holds.

Two self-dual additive codes over $\GF(9)$ are \emph{equivalent} if
the codewords of one can be mapped onto the codewords of the other by
a transformation that preserves the properties of the code, i.e.,
weight enumerator, additivity, and self-duality. It was shown by
Rains~\cite{nonbin1} that this group of transformations is $\symp_2(3)
\wr \sym(n)$, i.e., permutations of the coordinates combined with
operations from the \emph{symplectic group} $\symp_2(3)$ applied
independently to each coordinate.  Global conjugation of all
coordinates will also preserve the properties of the code, and codes
related by this operation are called \emph{weakly
  equivalent}~\cite{sdbook}.  In this paper, we classify codes up to
equivalence, i.e., we do not consider global conjugation.
Let an element
$a + b \omega \in \GF(9)$, be represented as $\binom{a}{b} \in \GF(3)^2$. We can then
premultiply this element by a $2 \times 2$ matrix. The group
\[
\symp_2(3) =
\left<\begin{pmatrix}1&1\\1&-1\end{pmatrix}, \begin{pmatrix}1&1\\0&1\end{pmatrix}\right>
\]
has order 24 and contains all $2 \times 2$ matrices with elements from
$\GF(3)$ and determinant one. The order of $\symp_2(3) \wr \sym(n)$ is
$24^n n!$, and hence this is the total number of maps that take a
self-dual additive code over $\GF(9)$ to an equivalent
code~\cite{nonbin1}.  By translating the action of $\symp_2(3)$ on
$\binom{a}{b}$ into operations on elements $c = a + b \omega \in
\GF(9)$, we find that the operations we can apply to all elements in a
coordinate of a code are $c \mapsto xc$ if $x^4 = 1$,
or $c \mapsto x\overline{c}$ if $x^4 = -1$, given $x \in \GF(9)$, and
$a + b \omega \mapsto a + yb + b \omega$, given $y \in \GF(3)$.

A transformation that maps $\mathcal{C}$ to itself is called an
\emph{automorphism} of $\mathcal{C}$.  All automorphisms of
$\mathcal{C}$ make up the \emph{automorphism group} of $\mathcal{C}$,
denoted $\Aut(\mathcal{C})$.  The number of distinct codes equivalent
to $\mathcal{C}$ is then given by $\frac{24^n
  n!}{\left|\Aut(\mathcal{C})\right|}$.  The \emph{equivalence class}
of $\mathcal{C}$ contains all distinct codes that are equivalent to
$\mathcal{C}$.  By adding the sizes of all equivalence classes of
codes of length~$n$, we find the total number of distinct codes of
length~$n$, denoted $T_n$. The number $T_n$ is also given by a
\emph{mass formula} which was described by Höhn~\cite{hohn} for
self-dual additive codes over $\GF(4)$ and is easily generalized to
$\GF(9)$:
\begin{equation}\label{eq:mass}
T_n = \prod_{i=1}^{n} (3^i+1) = \sum_{j=1}^{t_n} \frac{24^n
  n!}{\left|\Aut(\mathcal{C}_j)\right|},
\end{equation}
where $t_n$ is the number of equivalence classes of codes of
length~$n$, and $\mathcal{C}_j$ is a representative from each
equivalence class.
The smallest possible automorphism group, called the \emph{trivial
  automorphism group}, of a self-dual additive code over $\GF(9)$ is
$\{I,-I\}$, i.e., it consists of global multiplication of coordinates
by $1$ or $-1$.  By assuming that all codes of length $n$ have a
trivial automorphism group, we obtain from the mass formula a lower
bound on $t_n$, the total number of inequivalent codes.
\begin{equation}\label{eq:massbound}
t_n \ge \left\lceil\frac{2\prod_{i=1}^n (3^i+1)}{24^n n!}\right\rceil.
\end{equation}
Note that when $n$ is large, most codes have a trivial automorphism
group, so the tightness of the bound increases with $n$. As we will
see in Section~\ref{sec:conc}, for $n=10$, 80\% of all codes have a
trivial automorphism group, and the bound (\ref{eq:massbound})
underestimates $t_{10}$ by just 19\%.

Any \emph{linear} code over $\GF(9)$ that is self-dual with respect to
the \emph{Hermitian inner product}, $(\boldsymbol{u}, \boldsymbol{v})
= \boldsymbol{u} \cdot \overline{\boldsymbol{v}}$, is also a self-dual
additive code with respect to the Hermitian trace inner product.  The
class of Hermitian self-dual linear codes over $\GF(9)$ is also known
as $9^H$~\cite{sdbook}.  The operations that map a self-dual linear
code to an equivalent code are more restrictive than for additive
codes, since $\GF(9)$-linearity must now be preserved. Only coordinate
permutations and multiplication of single coordinates by $x \in
\GF(9)$ where $x^4 = 1$ is allowed.  It follows that only additive
codes that satisfy certain constraints can be equivalent to linear
codes.  Such constraints for codes over $\GF(4)$ were described by Van
den Nest~\cite{nestthesis} and by Glynn et~al.~\cite{glynnbook}. An
obvious constraint is that all coefficients of the weight enumerator,
except $A_0$, of a linear code must be divisible by 8, whereas for an
additive code they need only be divisible by 2.
To our knowledge, no complete classification of Hermitian self-dual
linear codes over $\GF(9)$ have appeared so far, but several authors
have studied this class of codes and suggested a number of
constructions~\cite{lin9gulliver,lin9kim,lin9grassl,lin9gaborit}.
Checking whether a self-dual additive code over $\GF(9)$ is equivalent
to a linear code is non-trivial, since there are $6^n$ coordinate
transformations in $\symp_2(3)^n$ that could transform a non-linear
code into a linear code. Our classification of self-dual additive
codes could be a useful starting point for also studying linear
codes, but this is left as a problem for future work.

Trace-Hermitian self-dual additive codes over $\GF(q)$ exist for
$q=m^2$, where $m$ is a prime power~\cite{sdbook}, and
the class of self-dual additive codes over $\GF(q)$ is
called $q^{H+}$. The first case, $4^{H+}$, has
been studied in detail, in particular since an application to
\emph{quantum error-correction} was discovered~\cite{calderbank}.  We
have previously classified self-dual additive codes over $\GF(4)$ up
to length 12~\cite{selfdualgf4}.  Self-dual linear codes over $\GF(4)$
have been classified up to length 16~\cite{gf4to16} by Conway, Pless,
and Sloane. This classification was recently extended to length
18~\cite{harada1} and 20~\cite{harada2} by Harada et al.  The next
class of self-dual additive codes, $9^{H+}$, has received less attention, 
although these codes have similar application in quantum
error-correction~\cite{nonbin1,nonbin2}, where they correspond to
\emph{ternary quantum codes}.  We have previously classified self-dual
additive codes over $\GF(9)$ up to length 8~\cite{nonbinary}, as well
as self-dual additive codes over $\GF(16)$ and $\GF(25)$ up to length 6. Another type of
self-dual code over $\GF(9)$ is known as $9^E$~\cite{sdbook} and is
self-dual with respect to the Euclidean inner product,
$(\boldsymbol{u}, \boldsymbol{v}) = \boldsymbol{u} \cdot
\boldsymbol{v}$. There is no additive variant of these codes, and this
family will not be considered in this paper. Again, some constructions
have been described~\cite{lin9grassl}, but no complete classifications of
Euclidean self-dual codes over $\GF(9)$ have been given.

In Section~\ref{sec:graph} we briefly review the connection between
trace-Hermitian self-dual additive codes and weighted graphs. An
algorithm for checking equivalence of self-dual additive codes over
$\GF(9)$, which is a generalization of a known algorithm for linear
codes~\cite{ostergard}, is described in Sections~\ref{sec:eqgraph}
and~\ref{sec:class}. Combining this algorithm with the weighted graph
representation, and some other optimizations, enables us to classify
all self-dual additive codes over $\GF(9)$ of length up to 10 in
Section~\ref{sec:classresult}.  In particular, all near-extremal codes
of length 9 and 10 are classified for the first time. We also find the
smallest codes with trivial automorphism group.  Using an extension
technique described in Section~\ref{sec:classopt}, we are then able to
classify all optimal codes of length 11 and 12.  
We finish with some concluding remarks in Section~\ref{sec:conc}.

\section{Codes and Weighted Graphs}\label{sec:graph}

An \emph{m-weighted graph} is a triple $G=(V,E,W)$, where $V$ is a set
of \emph{vertices}, $E \subseteq V \times V$ is a set of \emph{edges},
and $W$ is a set of weights from $\GF(m)$, such that each edge has an
associated non-zero weight. In an unweighted graph, which is simply
described by a pair $G=(V,E)$, we can consider all edges to have
weight one.  A graph with $n$ vertices can be represented by an $n
\times n$ \emph{adjacency matrix} $\Gamma$, where the element
$\Gamma_{i,j} = W({\{i,j\}})$ if $\{i,j\} \in E$, and $\Gamma_{i,j} =
0$ otherwise.  A \emph{loop-free} \emph{undirected} graph has a
symmetric adjacency matrix where all diagonal elements are 0. In a
\emph{directed} graph, edges are ordered pairs, and the adjacency
matrix is not necessarily symmetric. In a \emph{colored} graph, the
set of vertices is partitioned into disjoint subsets, where each
subset is assigned a different color.
Two graphs $G=(V,E)$ and $G'=(V,E')$ are \emph{isomorphic} if and only
if there exists a permutation $\pi$ of $V$ such that $\{u,v\} \in E
\iff \{\pi(u), \pi(v)\} \in E'$. For weighted graphs, we also require
that edge weights are preserved, i.e., $W(\{u,v\}) = W(\{\pi(u),
\pi(v)\})$. For a colored graph, we further require the permutation to
preserve the graph coloring, i.e., that all vertices are mapped to
vertices of the same color. The \emph{automorphism group} of a graph
is the set of vertex permutations that map the graph to itself. A
\emph{path} is a sequence of vertices, $(v_1,v_2,\ldots,v_i)$, such
that $\{v_1,v_2\}, \{v_2,v_3\},$ $\ldots, \{v_{i-1},v_{i}\} \in E$. A
graph is \emph{connected} if there is a path from any vertex to any
other vertex in the graph.

If an additive code over $\GF(9)$ has a generator matrix of the form
$C = \Gamma + \omega I$, where $I$ is the identity matrix, $\omega$ is
a primitive element of $\GF(9)$, and $\Gamma$ is the adjacency matrix
of a loop-free undirected 3-weighted graph, we say that the generator
matrix is in \emph{standard form}.
A generator matrix in standard form must generate a code that is
self-dual with respect to the Hermitian trace inner product, since it
has full rank over $\GF(3)$ and $C \overline{C}^{\text{T}} = \Gamma^2
+ \Gamma - I$ only contains entries from $\GF(3)$, and hence the
traces of all elements of $\omega^2 C \overline{C}^{\text{T}}$ will be zero.

In the context of quantum codes, it was shown by
Schlingemann~\cite{schlingemann2} and by Grassl, Klappenecker, and
Rötteler~\cite{grasslpaper} that every self-dual additive code is
equivalent to a code with a generator matrix in standard
form. Essentially, the same results was also shown by
Bouchet~\cite{bouchet} in the context of \emph{isotropic systems}.
The algorithm given in Fig.~\ref{fig:algcodetograph} can be used to
perform a mapping from a self-dual additive code to an equivalent code
in standard form.  Note that we can write the generator matrix $C = A
+ \omega B$ as an $n \times 2n$ matrix $(A\mid B)$ with elements from
$\GF(3)$. Steps 1 and 2 of the algorithm are used to obtain the
submatrices $A$ and $B$.  If $B$ now has full rank, we can simply
perform the basis change $B^{-1} (A \mid B) = (\Gamma \mid I)$ to
obtain the standard form.  Elements on the diagonal of $\Gamma$ can
then always be set to zero by operations $a + b \omega \mapsto a + yb
+ b \omega$, for $y \in \GF(3)$, corresponding to symplectic matrices
$\left(\begin{smallmatrix}1&y\\0&1\end{smallmatrix}\right)$. Hence
step 12 of the algorithm preserves code equivalence.  In the case
where $B$ has rank $k<n$, we can assume, after a basis change, that
the first $k$ rows and columns of $B$ form a $k \times k$ invertible
matrix. This is done in step 5, and the result is a permutation of the
coordinates of the code.  By the operation $c \mapsto
\omega\overline{c}$, for $c = a + b \omega$, corresponding to the
symplectic matrix
$\left(\begin{smallmatrix}0&-1\\1&0\end{smallmatrix}\right)$, we can
replace column $\boldsymbol{a_i}$ by $-\boldsymbol{b_i}$ and
$\boldsymbol{b_i}$ by $\boldsymbol{a_i}$. In this way, we ``swap'' the
$n-k$ last columns of $A$ and $B$ in steps 7 and 8. It has been shown
that it then follows from the self-duality of the code that the new
matrix $B$ must have full rank~\cite{nest,nonbinary}, and that the
matrix $\Gamma$ obtained in step 11 will always be symmetric.

\begin{figure}[t!]
\begin{algorithmic}[1]
  \medskip
  \REQUIRE $C$ generates a self-dual additive code over $\GF(9)$.
  \ENSURE $C'$ generates an equivalent code in standard form.
  \medskip
  \STATE $A \leftarrow \tr(\omega C)$
  \STATE $B \leftarrow \tr(\omega^2 C)$
  \STATE $k \leftarrow \rank(B)$
  \IF {$k < n$}
     \STATE Permute rows and columns of $B$ such that the first $k$ rows and columns form an invertible matrix. Apply the same permutation to the rows and columns of $A$.
     \FOR{$i=k+1$ to $n$}
        \STATE Swap columns $\boldsymbol{a_i}$ and $\boldsymbol{b_i}$
        \STATE $\boldsymbol{a_i} \leftarrow -\boldsymbol{a_i}$ 
     \ENDFOR
  \ENDIF
  \STATE $\Gamma \leftarrow B^{-1}A$
  \STATE Set all diagonal elements of $\Gamma$ to zero.
  \STATE $C' \leftarrow \Gamma + \omega I$
  \RETURN $C'$
  \medskip
\end{algorithmic}
 \caption{Algorithm for Mapping a Code to Standard Form}\label{fig:algcodetograph}
\end{figure}

\begin{figure}[t!]
 \centering
 \includegraphics[height=.25\linewidth]{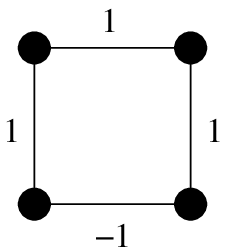}
 \caption{Graph Representation of the $(4,3^4,3)$ Code}\label{4d3graph}
\end{figure}

As an example, consider the $(4,3^4,3)$ code generated by $C$ which by
the described algorithm is transformed into the standard form
generator matrix $C'$, corresponding to the weighted graph depicted in
Fig.~\ref{4d3graph}:
\[
\setlength{\arraycolsep}{2.6pt}
C = 
\left(
\begin{array}{cccc}
1 & 0 & 1 & \omega^2  \\
\omega & 0 & \omega & \omega^3  \\
0 & 1 & \omega^2 & 1  \\
0 & \omega & \omega^3 & \omega
\end{array}
\right)
\,
C' = 
\left(
\begin{array}{cccc}
\omega & -1 & 0 & 1  \\
-1 & \omega & 1 & 0  \\
0 & 1 & \omega & 1  \\
1 & 0 & 1 & \omega
\end{array}
\right)
\]

It is known that two self-dual additive codes over $\GF(4)$ are
equivalent if and only if their corresponding graphs are related by a
sequence of graph operations called \emph{local complementations}
(LC)~\cite{bouchet,nest,glynnbook} and a permutation of the vertices.
We have previously used this fact to devise an algorithm to classify
all self-dual additive codes over $\GF(4)$ of length up to
12~\cite{selfdualgf4}.  The more general result that equivalence
classes of self-dual additive codes over $\GF(q=m^2)$ can be
represented as orbits of $m$-weighted graphs with respect to a
generalization of LC was later shown by Bahramgiri and
Beigi~\cite{newlc}. We used this to classify all self-dual additive
codes over $\GF(9)$, $\GF(16)$, and $\GF(25)$ up to lengths 8, 6, and
6, respectively~\cite{nonbinary}.  The main obstacle with this
approach is that the sizes of the LC orbits of weighted graphs quickly
get unmanageable as the number of vertices increase. We have therefore
devised a new method for checking code equivalence, which is described
in the next section. This algorithm uses a graph representation of
self-dual additive codes over $\GF(9)$ that is not related to the
representation described in this section, and does not require the
input to be in standard form. However, the weighted graph
representation will still be very useful for reducing our search
space.

\section{Equivalence Graphs}\label{sec:eqgraph}

\begin{figure}[!t]
\centering
\subfloat[]{\includegraphics[width=120pt]{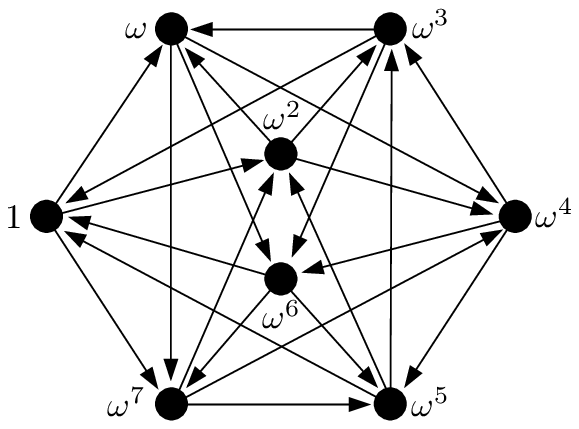}\label{subfig:sdadd}}
\subfloat[]{\includegraphics[width=120pt]{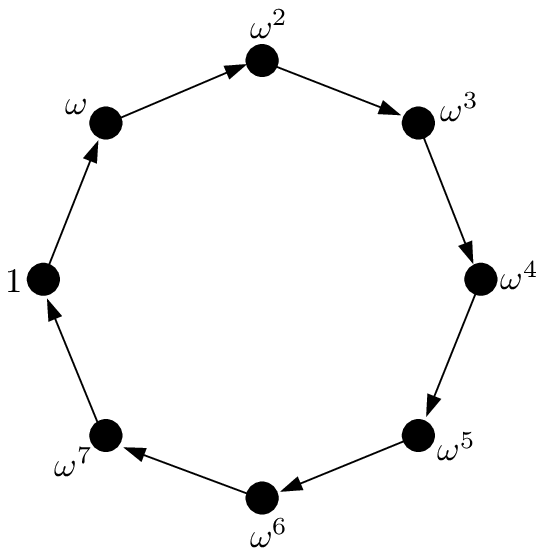}\label{subfig:linear}}\\
\hfil
\subfloat[]{\includegraphics[width=120pt]{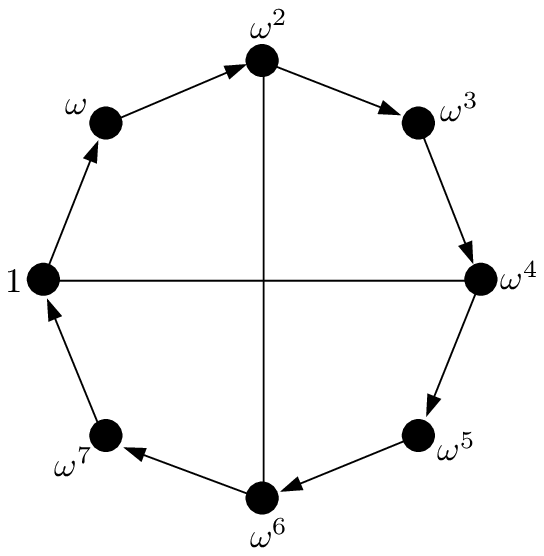}\label{subfig:hsd}}
\subfloat[]{\includegraphics[width=120pt]{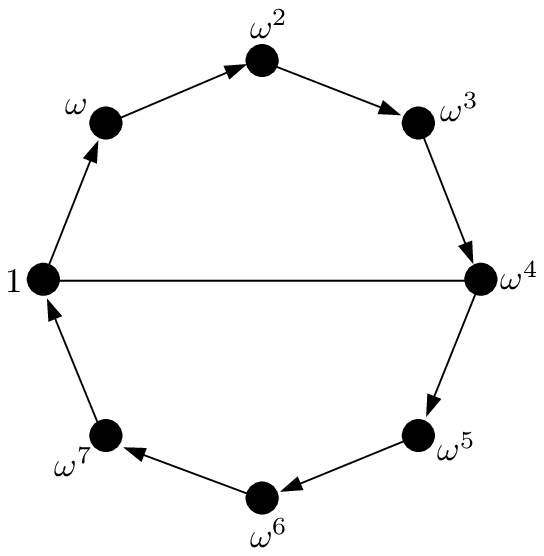}\label{subfig:esd}}
\hfil
\caption{Coordinate Graphs for Codes over GF(9): 
(a) Trace-Hermitian Self-Dual Additive 
(b) Linear 
(c) Hermitian Self-Dual Linear 
(d) Euclidean Self-Dual Linear}
\label{fig:cgraphs}
\end{figure}

To check whether two self-dual additive codes over $\GF(9)$ are
equivalent, we modify a well-known algorithm used for checking
equivalence of linear codes, described by
Östergård~\cite{ostergard}. The idea is to map a code to an
unweighted, directed, colored \emph{equivalence graph} such that the
automorphism groups of the code and the equivalence graph coincide.
An important component of the algorithm is to find a suitable
\emph{coordinate graph}. For self-dual additive codes over $\GF(9)$,
we need to construct a graph $G$ on eight vertices, labeled with the
non-zero elements of $\GF(9)$, whose automorphism group is
$\symp_2(3)$. This graph, shown in Fig.~\ref{fig:cgraphs}~\subref{subfig:sdadd}, was found
by adding directed edges $(\sigma 1, \sigma \omega)$ for all $\sigma
\in \symp_2(3)$. This ensures that $\symp_2(3) \subseteq \Aut(G)$. We
then verified that $\left|\Aut(G)\right| = 24$ which implies that
$\Aut(G) = \symp_2(3)$.

Fig.~\ref{fig:cgraphs} also shows examples of coordinate graphs for
some other families of codes over $\GF(9)$.  In the original algorithm
for checking equivalence of linear codes~\cite{ostergard}, the
coordinate graph shown in Fig.~\ref{fig:cgraphs}~\subref{subfig:linear} would be used. This
graph has an automorphism group of size eight, corresponding to the
fact that multiplication of a coordinate by any non-zero element of
$\GF(9)$ preserves linearity. The more restrictive coordinate graph
for Hermitian self-dual linear codes over $\GF(9)$ is shown in
Fig.~\ref{fig:cgraphs}~\subref{subfig:hsd}. This graph has an automorphism group of size
four, since only multiplication by $x \in \GF(9)$ where $x^4 = 1$ is
permitted in this case. Finally, Fig.~\ref{fig:cgraphs}~\subref{subfig:esd} shows a graph
with automorphism group of size two. This is the coordinate graph for
Euclidean self-dual linear codes over $\GF(9)$ where multiplication by
$\pm 1$ are the only permitted operations. Coordinate graphs of this
type were used by Harada and Östergård to classify Euclidean self-dual
codes over $\GF(5)$ up to length 16~\cite{ostergardharada1} and over
$\GF(7)$ up to length 12~\cite{ostergardharada2}.

\begin{figure}[!t]
\centering
\includegraphics[width=\columnwidth]{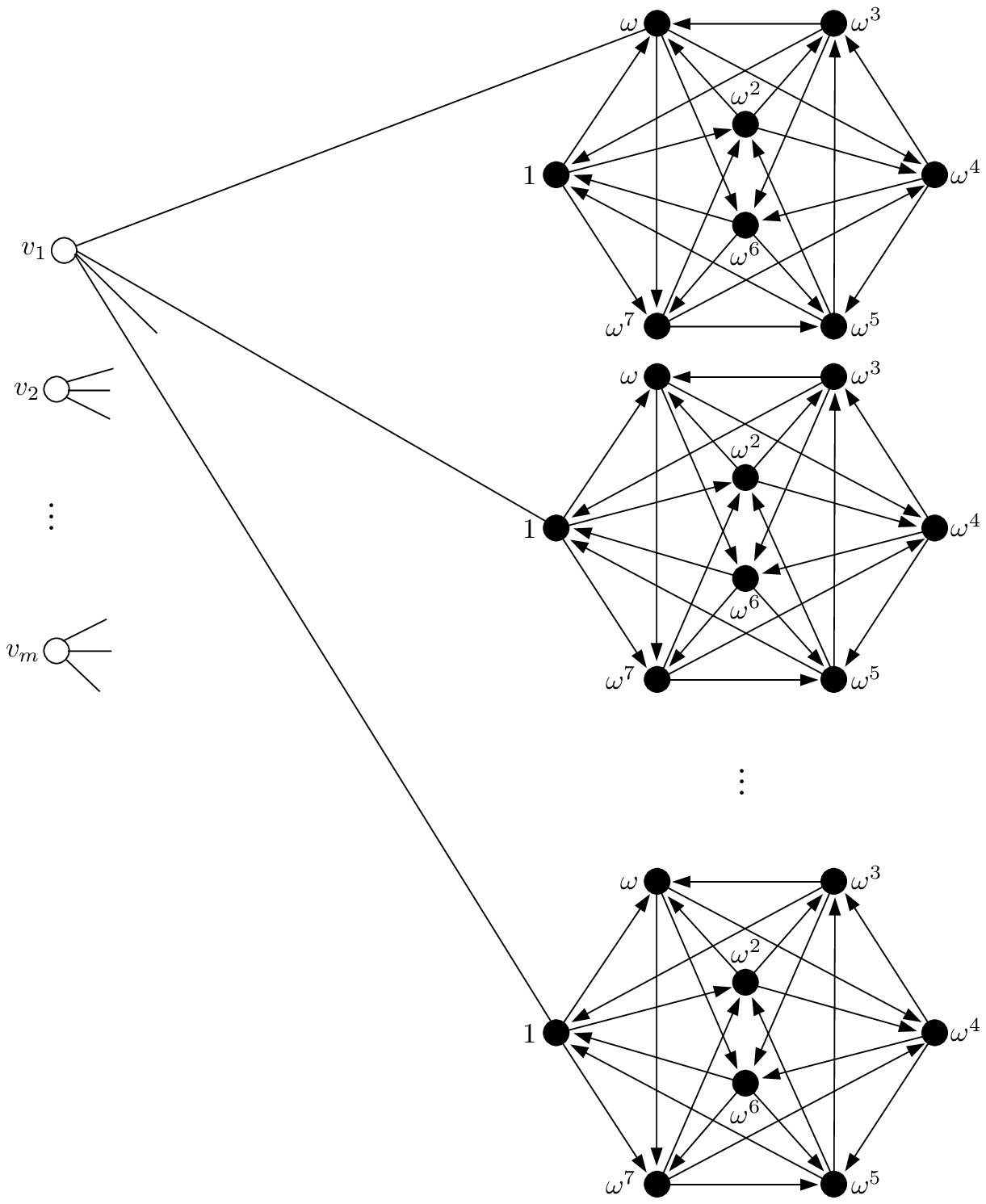}
\caption{Example of Equivalence Graph}
\label{fig:codegraph}
\end{figure}

To construct the equivalence graph of a code, we first add $n$ copies
of the coordinate graph, each copy representing one coordinate of the
code.  We then need a deterministic way to find a set of codewords
that generates the code. Taking all codewords would suffice, but the
following approach yields a smaller set and hence a more efficient
algorithm.  First, we check if the set of all codewords of minimum
weight $d$ generates the code. If it does not, we add all codewords of
weight $d+1$, then all codewords of weight $d+2$, etc, stopping once
we have a set that spans the code. For each codeword $c_i$ in the
resulting set, we add a \emph{codeword vertex} $v_i$ to the
equivalence graph. Let the codeword vertices have one color, and the
other vertices have a different color. Edges are added between $v_i$
and the coordinate graphs according to the non-zero coordinates of the
codeword $c_i$, e.g., if $c_i$ has $\omega$ in coordinate $j$, then
there is an edge between $v_i$ and the vertex labeled $\omega$ in the
$j$th coordinate graph. As an example, Fig.~\ref{fig:codegraph} shows
the case where $c_1 = (\omega 1 \cdots 1)$. The resulting equivalence
graph is finally \emph{canonized}, i.e., relabeled, but with coloring
preserved, using the \emph{nauty} software~\cite{nauty}. If two graphs
are isomorphic, their canonical representations are guaranteed to be
the same.

Applying a canonical permutation to the vertices of an equivalence
graph corresponds to permuting the coordinates of the corresponding
code, applying elements from $\symp_2(3)$ to each coordinate, and
sorting the codewords $c_i$ in some canonical order. If two codes are
equivalent, their canonical equivalence graphs will therefore be
identical. Furthermore, the automorphism group of a code is equivalent
to the automorphism group of its equivalence graph. This follows from
the fact that any automorphism of the equivalence graph must be one
out of $24^n n!$ possibilities, i.e., the $n!$ permutations of the $n$
coordinate subgraphs, and the 24 automorphisms from $\symp_2(3)$ of
each coordinate subgraph. No other automorphisms are possible. In
particular, permuting the codeword vertices will never be an
automorphism, since all codewords must be distinct.  Since it is
known~\cite{nonbin1} that coordinate permutations and $\symp_2(3)$
applied to the coordinates of a code preserve its weight enumerator,
additivity, and self-duality, this must also be true for any
automorphism of the equivalence graph.

\section{Classification Algorithm}\label{sec:class}

We have seen that every weighted graph corresponds to a self-dual
additive code, and that every self-dual additive code, up to
equivalence, has a standard form representation as a weighted
graph. It follows that we only need to consider 3-weighted graphs in
order to classify all self-dual additive codes over
$\GF(9)$. Permuting the vertices of a graph corresponds to permuting
coordinates of the associated code, which means that we only need to
consider these graphs up to isomorphism.  Moreover, we can restrict
our study to connected graphs, since a disconnected graph represents a
\emph{decomposable} code. A code is decomposable if it can be written
as the \emph{direct sum} of two smaller codes. For example, let
$\mathcal{C}$ be an $(n,3^n,d)$ code and $\mathcal{C}'$ an
$(n',3^{n'},d')$ code. The direct sum, $\mathcal{C} \oplus
\mathcal{C}' = \{u||v \mid u \in \mathcal{C}, v \in \mathcal{C}'\}$,
where $||$ means concatenation, is an $({n+n'},3^{n+n'},\min\{d,d'\})$
code. All decomposable codes of length~$n$ can be generated easily
once all indecomposable codes of length less than $n$ are known.

To classify codes of length $n$, we could take all non-isomorphic
connected 3-weighted graphs on $n$ vertices, map the corresponding
codes to equivalence graphs, and canonize these. All duplicates would
then be removed to obtain one representative code from each
equivalence class. However, a much smaller set of graphs is obtained
by taking all possible \emph{lengthenings}~\cite{gaborit} of all codes
of length $n-1$. A generator matrix in standard form can be lengthened
in $3^{n-1}-1$ ways by adding a vertex to the corresponding graph and
connecting it to all possible combinations of at least one of the old
vertices, using all possible combinations of edge weights. This
corresponds to adding a new non-zero row $\boldsymbol{r} \in \GF(3)^n$
and column $\boldsymbol{r}^\text{T}$ to the adjacency matrix, with
zero in the last coordinate.  Only half of the lengthenings need to be
considered, as adding the row $-\boldsymbol{r}$ is equivalent to
adding $\boldsymbol{r}$. (Since multiplying the last row and column in
the corresponding generator matrix by $-1$ would preserve code
equivalence.) We have previously shown~\cite{nonbinary}, using the
theory of local complementation of weighted graphs, that the set of
$i_{n-1}\frac{3^{n-1}-1}{2}$ codes obtained by lengthening one
representative from each of the $i_{n-1}$ equivalence classes of
indecomposable codes of length $n-1$ must contain at least one
representative from each equivalence class of the indecomposable codes
of length $n$.

Removing possible isomorphisms from the set of lengthened graphs,
using \emph{nauty}~\cite{nauty}, speeds up our classification
significantly. A set of non-isomorphic graphs that have already been
processed, as large as memory resources permit, can even be stored
between iterations, and new graphs can be checked for isomorphism
against this set. For each graph that is not excluded by such an
isomorphism check, the corresponding code must be mapped to an
equivalence graph, as described in Section~\ref{sec:eqgraph}. The
equivalence graph is canonized and compared against all previously
observed codes, which are stored in memory. Since the equivalence
graphs will be large, typically containing thousands of vertices for
$n=10$, we map the equivalence graph to a canonical generator matrix
by taking the first $n$ linearly independent codewords corresponding
to codeword vertices in their canonical ordering. This generator
matrix can further be mapped to a canonical standard form, as
described in Section~\ref{sec:graph}, which means that only
$\binom{n}{2}$ ternary symbols need to be stored for each code.  An
outline of the steps of our classification algorithm is listed in
Fig.~\ref{fig:algclassify}.

\begin{figure}[t!]
\begin{algorithmic}[1]
  \medskip
  \REQUIRE $\boldsymbol{C}_{n-1}$ contains one graph representation of each inequivalent indecomposable code of length $n-1$.
  \ENSURE $\boldsymbol{C}_{n}$ contains one graph representation of each inequivalent indecomposable code of length $n$.
  \medskip
  \STATE $\boldsymbol{C}_{n} \leftarrow \emptyset$
  \FORALL{$C \in \boldsymbol{C}_{n-1}$}
     \STATE $\boldsymbol{E} \leftarrow \frac{3^{n-1}-1}{2}$ lengthenings of $C$
     \STATE Remove isomorphisms from $\boldsymbol{E}$
     \FORALL{$E \in \boldsymbol{E}$}
        \STATE $d \leftarrow$ minimum distance of $E$
        \STATE $\boldsymbol{S} \leftarrow$ all codewords of weight $d$ from $E$
        \WHILE{$\boldsymbol{S}$ does not generate $E$}
           \STATE $d \leftarrow d+1$
           \STATE $\boldsymbol{S} \leftarrow \boldsymbol{S}\enspace \cup$ all codewords of weight $d$ from $E$
        \ENDWHILE
        \STATE $Q \leftarrow$ equivalence graph given by $\boldsymbol{S}$
        \STATE $Q' \leftarrow$ canonize $Q$
        \STATE $G \leftarrow$ graph representation of code given by $Q'$
        \IF{$G \not\in \boldsymbol{C}_{n}$}
           \STATE $\boldsymbol{C}_{n} \leftarrow \boldsymbol{C}_{n} \cup G$
        \ENDIF
     \ENDFOR
  \ENDFOR
  \RETURN $\boldsymbol{C}_n$
  \medskip
\end{algorithmic}
 \caption{Classification Algorithm}\label{fig:algclassify}
\end{figure}

Note that the special form of a generator matrix in standard form
makes it easy to find all codewords of low weight, which is necessary
to construct the equivalence graph. If $\mathcal{C}$ is generated by
$C = \Gamma + \omega I$, then any codeword formed by taking
$\GF(3)$-linear combinations of $i$ rows of $C$ must have weight at
least $i$. This means that we can find all codewords of weight $i$ by
only considering combinations of at most $i$ rows of $C$.

\section{Codes of Length 9 and 10}\label{sec:classresult}

Using the algorithm described in Section~\ref{sec:class}, we have
classified all self-dual additive codes over $\GF(9)$ of length up to
10. Table~\ref{tab:codes} gives the values of $i_n$, the number of
inequivalent indecomposable codes of length~$n$, and the values of
$t_n$, the total number of inequivalent codes of length~$n$.  Note
that the numbers $t_n$ are easily derived from the numbers $i_n$ by
using the \emph{Euler transform}~\cite{sloane2}:
\begin{eqnarray*}
c_n &=& \sum_{d|n} d i_d\\
t_1 &=& c_1\\
t_n &=& \frac{1}{n}\left( c_n + \sum_{k=1}^{n-1} c_k t_{n-k} \right).
\end{eqnarray*}
Tables~\ref{tab:distance} and~\ref{tab:distancetotal} list the numbers
of indecomposable codes and the total number of codes, respectively,
by length and minimum distance.  In Table~\ref{tab:wd}, we count the
number of distinct weight enumerators.  There are obviously too many
codes of length 9 and 10 to list all of them here, so an on-line
database containing one representative from each equivalence class has
been made available at \url{http://www.ii.uib.no/~larsed/nonbinary/}.

\begin{table}[!t]
\renewcommand{\arraystretch}{1.3}
\setlength{\extrarowheight}{1pt}
\caption{Number of Indecomposable ($i_n$) and Total Number ($t_n$) of Self-Dual Additive Codes over GF(9) 
of Length~$n$}%
\label{tab:codes}
\centering
\begin{tabular}{crrrrrrrrrr}
\hline
$n$   & 1 & 2 & 3 & 4 &  5 &  6 &   7 &   8 &      9 &         10 \\
\hline
$i_n$ & 1 & 1 & 1 & 3 &  5 & 21 &  73 & 659 & 17\,589 &  2\,803\,404 \\
$t_n$ & 1 & 2 & 3 & 7 & 13 & 39 & 121 & 817 & 18\,525 &  2\,822\,779 \\
\hline
\end{tabular}
\end{table}

\begin{table*}
\renewcommand{\arraystretch}{1.3}
\setlength{\extrarowheight}{1pt}
\caption{Number of Indecomposable Self-Dual Additive Codes over GF(9) of Length~$n$ and Minimum Distance~$d$}%
\label{tab:distance}
\centering
\begin{tabular}{crrrrrrrrrrr}
\hline
$d \backslash n$ & 
2 & 3 & 4 & 5 & 6 & 7 & 8 & 9 & 10 & 11 & 12 \\
\hline
2     & 1 & 1 & 2 &  4 & 15 & 51 & 388 &    6240 &    418\,088 & ?            & ? \\
3     &   &   & 1 &  1 &  5 & 20 & 194 &    6975 &    893\,422 & ?            & ? \\
4     &   &   &   &    &  1 &  2 &  77 &    4370 & 1\,487\,316 & ?            & ? \\
5     &   &   &   &    &    &    &     &       4 &        4577 & 56\,005\,876 & ? \\
6     &   &   &   &    &    &    &     &         &           1 &              & 6493 \\
\hline
All   & 1 & 1 & 3 &  5 & 21 & 73 & 659 & 17\,589 & 2\,803\,404 & ?            & ? \\
\hline
\end{tabular}
\end{table*}

\begin{table*}
\renewcommand{\arraystretch}{1.3}
\setlength{\extrarowheight}{1pt}
\caption{Total Number of Self-Dual Additive Codes over GF(9) of Length~$n$ and Minimum Distance~$d$}%
\label{tab:distancetotal}
\centering
\begin{tabular}{crrrrrrrrrrrr}
\hline
$d \backslash n$ & 
1 & 2 & 3 & 4 & 5 & 6 & 7 & 8 & 9 & 10 & 11 & 12 \\
\hline
1     & 1 & 1 & 2 & 3 &  7 & 13 &  39 & 121 &     817 &     18\,525 & 2\,822\,779  & ? \\
2     &   & 1 & 1 & 3 &  5 & 20 &  60 & 424 &    6358 &    418\,931 & ?            & ? \\
3     &   &   &   & 1 &  1 &  5 &  20 & 195 &    6976 &    893\,429 & ?            & ? \\
4     &   &   &   &   &    &  1 &   2 &  77 &    4370 & 1\,487\,316 & ?            & ? \\
5     &   &   &   &   &    &    &     &     &       4 &        4577 & 56\,005\,876 & ? \\
6     &   &   &   &   &    &    &     &     &         &           1 &              & 6493 \\
\hline 
All   & 1 & 2 & 3 & 7 & 13 & 39 & 121 & 817 & 18\,525 & 2\,822\,779 & $> 2^{30}$  & $> 2^{41}$ \\
\hline
\end{tabular}
\end{table*}

\begin{table*}
\renewcommand{\arraystretch}{1.3}
\setlength{\extrarowheight}{1pt}
\caption{Number of Distinct Weight Enumerators of Indecomposable Codes of Length~$n$ and Minimum Distance~$d$}%
\label{tab:wd}
\centering
\begin{tabular}{crrrrrrrrrrrr}
\hline
$d \backslash n$ & 
2 & 3 & 4 & 5 & 6 & 7 & 8 & 9 & 10 & 11 & 12 \\
\hline
2     & 1 & 1 & 2 &  4 & 14 & 42 & 202 &  1021 &  8396 &  ? & ? \\
3     &   &   & 1 &  1 &  3 &  9 &  33 &   170 &  1133 &  ? & ? \\
4     &   &   &   &    &  1 &  1 &   9 &    25 &   345 &  ? & ? \\
5     &   &   &   &    &    &    &     &     1 &    10 & 48 & ? \\
6     &   &   &   &    &    &    &     &       &     1 &    & 27 \\
\hline
All   & 1 & 1 & 3 &  5 & 18 & 52 & 244 &  1217 &  9885 &  ? & ? \\
\hline
\end{tabular}
\end{table*}

Generator matrices for all the extremal codes of length 9 and 10 were
given in~\cite{nonbinary}.  Our classification confirms that there are
four extremal $(9,3^9,5)$ codes, all with weight enumerator $W(y) = 1
+ 252 y^5 + 1176 y^6 + 3672 y^7 + 7794 y^8 + 6788 y^9$, and with
automorphism groups of size 72, 108, 108, and 432, and that there is a
unique extremal $(10,3^{10},6)$ code with weight enumerator $W(y) = 1
+ 1680 y^6 + 2880 y^7 + 14040 y^8 + 22160 y^9 + 18288 y^{10}$ and
automorphism group of size 2880.  The classification of near-extremal
codes of length 9 and 10 is new.  We find that there are 4370
near-extremal $(9,3^9,4)$ codes with 25 distinct weight enumerators
and 13 different values for $\left|\Aut(\mathcal{C})\right|$.  The
weight enumerators that exist are given by $W_{9,\alpha}(y) = 1 +
(4+2\alpha)y^4 + (244-4\alpha)y^5 + (1168-4\alpha)y^6 +
(3704+16\alpha)y^7 + (7766-14\alpha)y^8 + (6796+4\alpha)y^9$ for all
integer values $0 \le \alpha \le 24$.  Table~\ref{tab:nearext9} gives
the number of $(9,3^9,4)$ codes for each possible weight enumerator
and automorphism group size.  To highlight a few codes with extreme
properties, we list generator matrices for the code with automorphism
group of maximal size (288) and one of the codes with weight
enumerator $W_{9,0}(y)$, i.e., with the minimal number of minimum
weight codewords. For the latter case, we choose the unique code with
maximum number of automorphisms (16). In the following, ``-'' denotes
$-1$ in generator matrices:

\begin{table*}
\renewcommand{\arraystretch}{1.3}
\setlength{\extrarowheight}{1pt}
\caption{Number of $(9,3^9,4)$ Codes with Weight Enumerator~$W_{9,\alpha}(y)$ and
  $\left|\Aut(\mathcal{C})\right| = \beta$}%
\label{tab:nearext9}
\centering
\begin{tabular}{crrrrrrrrrrrrrr}
\hline
$\alpha \backslash \beta$ & 
2 & 4 & 6 & 8 & 12 & 16 & 24 & 32 & 36 & 48 & 72 & 144 & 288 & All \\
\hline
0     &     &   3 &   &  1 &    &  1 &   &   &  &  &  &  &  &   5 \\
1     &   2 &   1 & 2 &    &    &    &   &   &  &  &  &  &  &   5 \\
2     &  15 &  21 &   &  4 &    &    &   &   &  &  &  &  &  &  40 \\
3     &  15 &  13 & 1 &  3 &  2 &    & 1 &   &  &  &  &  &  &  35 \\
4     & 125 &  52 &   &    & 12 &  2 & 2 &   &  & 2&  &  &  & 195 \\
5     &  85 &   8 &   &    &    &    &   &   &  &  &  &  &  &  93 \\
6     & 338 &  93 &   & 11 &  2 &  1 &   &   &  &  &  &  &  & 445 \\
7     & 165 &  53 & 2 &  9 &  2 &  2 & 2 &   &  &  & 1&  &  & 236 \\
8     & 561 & 150 &   & 11 &    &    &   & 1 &  &  &  &  &  & 723 \\
9     & 173 &  20 & 6 &  7 &    &    & 1 &   &  &  &  &  &  & 207 \\
10    & 522 & 154 & 4 &  7 &  7 &    & 2 &   &  & 2&  &  &  & 698 \\
11    & 157 &  53 &   & 15 &    &    &   &   &  &  &  &  &  & 225 \\
12    & 356 & 143 & 2 &  4 &  3 &  2 &   &   &  &  &  &  &  & 510 \\
13    & 119 &  25 & 2 &  6 &    &    &   &   &  &  &  &  &  & 152 \\
14    & 229 & 114 &   & 11 &    &  2 &   &   &  &  &  &  &  & 356 \\
15    &  42 &  28 & 1 & 16 &  1 &  2 &   &   &  &  &  &  &  &  90 \\
16    &  96 &  62 &   &  8 &  2 &  3 & 2 &   & 4&  &  & 2& 1& 180 \\
17    &  15 &   9 &   &  6 &    &    &   &   &  &  &  &  &  &  30 \\
18    &  23 &  33 &   &  2 &  1 &  6 &   & 2 &  &  &  &  &  &  67 \\
19    &   9 &   4 &   &  6 &    &  2 & 2 &   &  &  &  &  &  &  23 \\
20    &   8 &  23 &   &  6 &    &    &   & 2 &  &  &  &  &  &  39 \\
21    &     &   2 & 2 &    &    &    & 1 &   &  &  &  &  &  &   5 \\
22    &   1 &   3 &   &    &    &  2 &   & 1 &  &  &  &  &  &   7 \\
23    &     &     &   &  1 &    &    &   &   &  &  &  &  &  &   1 \\
24    &     &     &   &    &    &    &   & 3 &  &  &  &  &  &   3 \\
\hline
All   & 3056 & 1067 & 22 & 134 & 32 & 25 & 13 & 9 & 4 & 4 & 1 & 2 & 1 & 4370 \\
\hline
\end{tabular}
\end{table*}

\[\small
\setlength{\arraycolsep}{1.5pt}
C^{n=9}_{\left|\Aut\right|=288} =
\left(
\begin{array}{ccccccccc}
\omega&$-$&1&1&$-$&$-$&1&$-$&1\\
$-$&\omega&1&1&1&0&1&$-$&0\\
1&1&\omega&1&1&0&0&1&$-$\\
1&1&1&\omega&1&1&$-$&0&0\\
$-$&1&1&1&\omega&$-$&0&0&1\\
$-$&0&0&1&$-$&\omega&0&1&0\\
1&1&0&$-$&0&0&\omega&0&1\\
$-$&$-$&1&0&0&1&0&\omega&0\\
1&0&$-$&0&1&0&1&0&\omega
\end{array}
\right)
\]

\[\small
\setlength{\arraycolsep}{1.5pt}
C^{n=9}_{\alpha=0,\, \left|\Aut\right|=16} =
\left(
\begin{array}{ccccccccc}
\omega&1&1&1&1&$-$&$-$&$-$&$-$\\
1&\omega&1&$-$&$-$&1&$-$&1&0\\
1&1&\omega&$-$&$-$&$-$&1&0&1\\
1&$-$&$-$&\omega&1&0&1&$-$&1\\
1&$-$&$-$&1&\omega&1&0&1&$-$\\
$-$&1&$-$&0&1&\omega&$-$&$-$&1\\
$-$&$-$&1&1&0&$-$&\omega&1&$-$\\
$-$&1&0&$-$&1&$-$&1&\omega&$-$\\
$-$&0&1&1&$-$&1&$-$&$-$&\omega
\end{array}
\right)
\]

We find that there are 4577
near-extremal $(10,3^{10},5)$ codes with 10 distinct weight
enumerators and 20 different values for
$\left|\Aut(\mathcal{C})\right|$. The weight enumerators that exist
are given by $W_{10,\alpha}(y) = 1 + (44+4\alpha)y^5 +
(1460-20\alpha)y^6 + (3320+40\alpha)y^7 + (13600-40\alpha)y^8 +
(22380+20\alpha)y^9 + (18244-4\alpha)y^{10}$ for integer values
$\alpha \in \{0, 9, 12, 13, 16, 18, 21, 22, 24, 25\}$.
Table~\ref{tab:nearext10} gives the number of $(10,3^{10},5)$ codes
for each possible weight enumerator and automorphism group size.
We give generator matrices for the unique codes with automorphism groups of size 2880
and 288, as well as the unique code with weight enumerator $W_{10,0}(y)$:

\begin{table*}
\renewcommand{\arraystretch}{1.3}
\setlength{\extrarowheight}{1pt}
\caption{Number of $(10,3^{10},5)$ Codes with Weight Enumerator~$W_{10,\alpha}(y)$ and
  $\left|\Aut(\mathcal{C})\right| = \beta$}%
\label{tab:nearext10}
\centering
\begin{tabular}{crrrrrrrrrrrrrrrrrrrrr}
\hline
$\alpha \backslash \beta$ & 
      2  & 4 & 6 & 8 & 10 & 12 & 16 & 20 & 24 & 32 & 36 & 40 & 48 & 64 & 72 & 144 & 192 & 240 & 288 & 2880 & All  \\
\hline
0   &    &   &   &   &    &    &    &    &    &    &    &    &    &    &    &     &     &    1&     &      &    1 \\
9   &    &   &   &   &    &    &   1&    &    &    &    &   1&    &    &    &     &    1&    1&     &      &    4 \\
12  &    &   &   &   &    &    &   2&    &    &    &    &    &    &    &    &     &     &     &     &      &    2 \\
13  &    &  3&   &   &    &    &   3&    &    &   1&    &    &    &    &    &     &     &     &     &      &    7 \\
16  &  10&  5&   &  2&    &    &    &    &    &    &   1&    &    &    &    &     &     &     &     &      &   18 \\
18  &  30& 24&  4&  4&    &   8&   4&    &    &    &    &    &    &    &    &     &     &     &     &      &   74 \\
21  & 190& 77&  2& 20&    &   2&   2&    &   4&    &    &    &    &   3&    &     &     &     &     &      &  300 \\
22  & 467& 72&   &  4&    &    &   1&    &    &    &    &    &    &    &    &     &     &     &     &      &  544 \\
24  &2321&172&  4&  4&   1&   5&   1&    &    &    &    &   1&   1&    &    &     &     &     &     &      & 2510 \\
25  & 777&247& 12& 39&    &  14&  10&   3&   2&   1&   2&   2&    &   4&   1&    1&     &     &    1&     1& 1117 \\
\hline
All & 3795 & 600 & 22 & 73 & 1 & 29 & 24 & 3 & 6 & 2 & 3 & 4 &  1 &  7 &  1 &   1 &   1 &   2 &   1 &    1 & 4577 \\
\hline
\end{tabular}
\end{table*}

\[\small
\setlength{\arraycolsep}{1.5pt}
C^{n=10}_{\left|\Aut\right|=2880} =
\left(
\begin{array}{cccccccccc}
\omega&1&1&1&1&0&$-$&$-$&1&1\\
1&\omega&1&1&1&1&$-$&1&$-$&0\\
1&1&\omega&1&1&1&1&$-$&0&$-$\\
1&1&1&\omega&1&$-$&0&1&1&$-$\\
1&1&1&1&\omega&$-$&1&0&$-$&1\\
0&1&1&$-$&$-$&\omega&1&1&1&1\\
$-$&$-$&1&0&1&1&\omega&1&1&1\\
$-$&1&$-$&1&0&1&1&\omega&1&1\\
1&$-$&0&1&$-$&1&1&1&\omega&1\\
1&0&$-$&$-$&1&1&1&1&1&\omega
\end{array}
\right)
\]

\[\small
\setlength{\arraycolsep}{1.5pt}
C^{n=10}_{\left|\Aut\right|=288} =
\left(
\begin{array}{cccccccccc}
\omega&0&0&$-$&1&0&0&0&1&1\\
0&\omega&$-$&$-$&1&0&1&0&$-$&1\\
0&$-$&\omega&$-$&1&0&0&1&1&$-$\\
$-$&$-$&$-$&\omega&1&0&$-$&$-$&1&1\\
1&1&1&1&\omega&0&0&0&0&0\\
0&0&0&0&0&\omega&1&$-$&$-$&1\\
0&1&0&$-$&0&1&\omega&$-$&1&0\\
0&0&1&$-$&0&$-$&$-$&\omega&0&1\\
1&$-$&1&1&0&$-$&1&0&\omega&$-$\\
1&1&$-$&1&0&1&0&1&$-$&\omega
\end{array}
\right)
\]

\[\small
\setlength{\arraycolsep}{1.5pt}
C^{n=10}_{\alpha=0} =
\left(
\begin{array}{cccccccccc}
\omega&1&0&$-$&$-$&1&0&$-$&0&0\\
1&\omega&1&1&0&1&0&$-$&1&$-$\\
0&1&\omega&0&1&1&0&1&1&0\\
$-$&1&0&\omega&0&$-$&1&$-$&$-$&$-$\\
$-$&0&1&0&\omega&1&$-$&$-$&$-$&$-$\\
1&1&1&$-$&1&\omega&1&$-$&1&0\\
0&0&0&1&$-$&1&\omega&1&$-$&0\\
$-$&$-$&1&$-$&$-$&$-$&1&\omega&1&1\\
0&1&1&$-$&$-$&1&$-$&1&\omega&1\\
0&$-$&0&$-$&$-$&0&0&1&1&\omega
\end{array}
\right)
\]

That our classification of all codes up to length 10 is correct has
been verified by the mass formula (\ref{eq:mass}). This required us to
also calculate the sizes of the automorphism groups of all
decomposable codes, which was simplified by the observation that for a
code $\mathcal{C} = k_1\mathcal{C}_1 \oplus \cdots \oplus
k_m\mathcal{C}_m$, where $k_j\mathcal{C}_j = \bigoplus_{i=1}^{k_j}
\mathcal{C}_j$, $\left|\Aut(\mathcal{C})\right| = \prod_{i=1}^m
k_i!\left|\Aut(\mathcal{C}_i)\right|^{k_i}$.

Table~\ref{tab:trivial} gives the numbers of codes with trivial
automorphism group by length and minimum distance.  We find that the
smallest codes with trivial automorphism group are 35 codes of length
8. (Note that automorphism group sizes were not calculated in the
previous classification of codes of length 8~\cite{nonbinary}.)  We
give the generator matrix for one $(8,3^8,4)$ code with trivial
automorphism group. Generator matrices for the other codes can be
obtained from \url{http://www.ii.uib.no/~larsed/nonbinary/}.

\begin{table}[!t]
\renewcommand{\arraystretch}{1.3}
\setlength{\extrarowheight}{1pt}
\caption{Number of Codes of Length~$n$ and Minimum Distance~$d$ with Trivial Automorphism Group}%
\label{tab:trivial}
\centering
\begin{tabular}{crrrrrr}
\hline
$d \backslash n$ & 
$\le 7$ & 8 & 9 & 10 & 11 & 12 \\
\hline
 $\le2$& 0 &  0 &     0 &         0   &            0 & 0 \\
 3     & 0 & 32 &  4518 &   832\,878  &            ? & ? \\
 4     & 0 &  3 &  3056 & 1\,419\,861 &            ? & ? \\
 5     &   &    &     0 &     3795    & 55\,865\,753 & ? \\
 6     &   &    &       &        0    &              & 3445 \\
\hline
 All   & 0 & 35 &  7574 & 2\,256\,534 &            ? & ? \\
\hline
\end{tabular}
\end{table}

\[\small
\setlength{\arraycolsep}{1.5pt}
C^{n=8}_{\left|\Aut\right|=2} = 
\left(
\begin{array}{cccccccc}
\omega&0&0&1&1&$-$&$-$&$-$\\
0&\omega&$-$&0&1&1&1&$-$\\
0&$-$&\omega&$-$&1&0&$-$&0\\
1&0&$-$&\omega&0&$-$&0&1\\
1&1&1&0&\omega&1&0&0\\
$-$&1&0&$-$&1&\omega&0&0\\
$-$&1&$-$&0&0&0&\omega&$-$\\
$-$&$-$&0&1&0&0&$-$&\omega
\end{array}
\right)
\]

We observe that codes with minimum distance $d \le 2$ always have
nontrivial automorphisms, and this can be proved as follows. For
$d=1$, we can assume that the first row of a standard form generator
matrix is $(\omega 0 \cdots 0)$. Then
$\left(\begin{smallmatrix}1&0\\1&1\end{smallmatrix}\right)$ applied to
the first coordinate of the code is an automorphism of
order~3. Multiplying the first coordinate by $-1$ has the same effect
as multiplying the first row by $-1$ and is therefore an automorphism
of order~2.  Including the trivial automorphism, we have that
$\left|\Aut\right| \ge 12$. There are codes of length 9 with $d=1$ and
$\left|\Aut\right| = 12$ which shows that the bound is tight. For
$d=2$, we can assume that the first row of a standard form generator
matrix is $(\omega 1 0 \cdots 0)$. Then
$\left(\begin{smallmatrix}1&0\\1&1\end{smallmatrix}\right)$ applied to
the first coordinate and
$\left(\begin{smallmatrix}1&1\\0&1\end{smallmatrix}\right)$ applied to
the second coordinate of the code has the same effect as adding the
first row of the generator matrix to the second row, and is hence an
automorphism of order~3.  Swapping the first two coordinates is an
automorphism of order~2, since it has the same effect as the following
procedure: Add the first row to itself, then add the first row to each
row $i > 2$ where the value in position $i$ of the second column is 1,
and add twice the first row to each row $i > 2$ where the value in
position $i$ of the second column is 2. Finally apply
$\left(\begin{smallmatrix}0&2\\1&0\end{smallmatrix}\right)$ to the
first column, and
$\left(\begin{smallmatrix}0&1\\2&0\end{smallmatrix}\right)$ to the
second column.  Again, including the trivial automorphism we get the
bound $\left|\Aut\right| \ge 12$, and the existence of codes of length
8 with $d=2$ and $\left|\Aut\right| = 12$ proves that the bound is
tight.

\section{Optimal Codes of Length 11 and 12}\label{sec:classopt}

When we lengthen an $(n,3^n,d)$ code, as described in
Section~\ref{sec:class}, we always obtain an $(n+1,3^{n+1},d')$ code
where $d' \le d+1$~\cite{gaborit}. It follows that given a
classification of all codes of length $n$ and minimum distance at
least $d$, we can classify all codes of length $n+1$ and minimum
distance at least $d+1$. There are no $(11,3^{11},6)$ codes, but by
lengthening the 1\,491\,894 $(10,3^{10},d)$ codes for $d \ge 4$, we
are able to obtain all optimal $(11,3^{11},5)$ codes.  To quickly
exclude codes with $d < 5$, we checked the minimum distance of each
lengthened code before checking for code equivalence in this search.

We find that there are 56\,005\,876 optimal $(11,3^{11},5)$ codes with
48 distinct weight enumerators and 24 different values for
$\left|\Aut(\mathcal{C})\right|$. The weight enumerators that exist
are given by $W_{11,\alpha}(y) = 1 + (12+2\alpha)y^5 +
(888-6\alpha)y^6 + 3960y^7 + (14970+20\alpha)y^8 + (42500-30\alpha)y^9
+ (66240+18\alpha)y^{10} + (48756-4\alpha)y^{11}$ for all integer
values $6 \le \alpha \le 50$ as well as $\alpha \in \{0, 54, 60\}$.
Observe that the number of codewords of weight 7 is constant for all
codes. Table~\ref{tab:opt11} gives the number of $(11,3^{11},5)$ codes
for each possible weight enumerator and automorphism group size.
We give generator matrices for the unique codes with automorphism group
of size 47\,520 and 1440, as well as the unique code with weight enumerator
$W_{11,0}(y)$:

\begin{table*}
\renewcommand{\arraystretch}{1.3}
\setlength{\tabcolsep}{3.25pt}
\setlength{\extrarowheight}{1pt}
\caption{Number of $(11,3^{11},5)$ Codes with Weight Enumerator~$W_{11,\alpha}(y)$ and
  $\left|\Aut(\mathcal{C})\right| = \beta$}%
\label{tab:opt11}
\centering
\begin{tabular}{crrrrrrrrrrrrrrrrrrrrrrrrr}
\hline
$\alpha \backslash \beta$ & 
              2&      4& 6& 8&10&12&16&18&20&24&32&36&40&44&48&72&108&120&144&288&360&432&1440&47\,520& All  \\
\hline
 0 &           &       &  &  &  &  &  & 1&  &  &  &  &  &  &  &  &   &   &   &   &   &   &    &  &           1 \\
 6 &           &       & 4&  &  &  &  & 1&  &  &  & 3&  &  &  & 2&   &   &   &   &   &   &    &  &          10 \\
 7 &           &      4&  & 1&  &  &  &  &  & 2&  &  &  &  &  &  &   &   &   &   &   &   &    &  &           7 \\
 8 &         10&      5&  &  &  &  &  &  &  &  &  &  &  &  &  &  &   &   &   &   &   &   &    &  &          15 \\
 9 &         36&     22& 4& 4&  & 4&  &  &  &  &  & 2&  &  &  &  &   &   &   &   &   &   &    &  &          72 \\
10 &         35&     16& 2& 2&  & 2&  &  &  &  &  &  &  &  &  &  &   &   &   &   &   &   &    &  &          57 \\
11 &        286&     62&  & 2&  &  &  &  &  &  &  &  &  &  &  &  &   &   &   &   &   &   &    &  &         350 \\
12 &        217&     37&15& 7&  & 2& 1&  &  &  &  &  &  &  &  &  &   &   &   &   &   &   &    &  &         279 \\
13 &       1515&    170& 6& 8&  & 4&  &  &  & 2&  &  &  &  &  &  &   &   &   &   &   &   &    &  &        1705 \\
14 &       1140&    139&  & 4&  &  &  &  &  &  &  &  &  &  &  &  &   &   &   &   &   &   &    &  &        1283 \\
15 &       7412&    414&10&20&  &14&  &  & 1&  &  & 8&  &  &  & 5&   &   &   &   &   &   &    &  &        7884 \\
16 &       5234&    192& 4&13&  &  & 2&  &  &  &  &  &  &  &  &  &   &   &   &   &   &   &    &  &        5445 \\
17 &    30\,825&    906&  &28&  &  &  &  &  &  &  &  &  &  &  &  &   &   &   &   &   &   &    &  &     31\,759 \\
18 &    19\,468&    623&17&14&  &12&  &  &  &  &  &  &  &  &  &  &   &   &   &   &   &   &    &  &     20\,134 \\
19 &   108\,109&   1606&26&24&  & 8&  &  &  &  &  &  &  &  &  &  &   &   &   &   &   &   &    &  &    109\,773 \\
20 &    62\,364&    641&  & 7&  & 5&  &  &  &  &  &  &  &  &  &  &   &   &   &   &   &   &    &  &     63\,017 \\
21 &   314\,156&   2701&16&42&  &27&  &  &  & 2&  &  &  &  &  &  &   &   &   &   &   &   &    &  &    316\,944 \\
22 &   169\,270&   1928& 4&30&  &10&  &  &  &  &  &  &  &  &  &  &   &   &   &   &   &   &    &  &    171\,242 \\
23 &   780\,271&   4123& 5&40&  &  &  &  &  &  &  &  &  &  &  &  &   &   &   &   &   &   &    &  &    784\,439 \\
24 &   385\,400&   1508&38&32& 1& 8& 2&  & 1&  & 2& 3&  &  &  & 1&   &   &  1&  1&   &   &   1&  &    386\,999 \\
25 &1\,649\,942&   5666&42&33&  & 2&  &  &  &  &  &  &  &  &  &  &   &   &   &   &   &   &    &  & 1\,655\,685 \\
26 &   754\,931&   4249&  &44&  &10& 3&  &  & 4&  &  &  &  &  &  &   &   &   &   &   &   &    &  &    759\,241 \\
27 &2\,990\,527&   7882&61&36&  &44&  & 2&  &  &  & 2&  & 6&  &  &   &   &   &   &   &   &    &  & 2\,998\,560 \\
28 &1\,266\,193&   2610&20&19&  & 6&  &  &  &  & 1&  &  &  & 2&  &   &   &   &   &   &   &    &  & 1\,268\,851 \\
29 &4\,671\,482&   9256&18&36&  &  &  &  &  &  &  &  & 2&  &  &  &   &   &   &   &   &   &    &  & 4\,680\,794 \\
30 &1\,832\,724&   6641&41&50&  &20&  &  &  &  &  &  &  &  &  &  &   &   &   &   &   &   &    &  & 1\,839\,476 \\
31 &6\,241\,827&10\,336&98&39&  &10&  &  &  &  &  &  &  &  &  &  &   &   &   &   &   &   &    &  & 6\,252\,310 \\
32 &2\,266\,449&   3048&  &45&  &  & 6&  &  &  &  &  &  &  &  &  &   &   &   &   &   &   &    &  & 2\,269\,548 \\
33 &7\,110\,043&10\,986&89&27&  &47&  &  &  &  &  & 7&  &  &  & 4&  2&   &   &   &   &   &    &  & 7\,121\,205 \\
34 &2\,377\,017&   7970&44&66&  & 4& 6&  &  &  &  &  &  &  &  &  &   &   &   &   &   &   &    &  & 2\,385\,107 \\
35 &6\,821\,413&10\,684& 6&22&  &  &  &  &  &  &  &  &  &  &  &  &   &   &   &   &   &   &    &  & 6\,832\,125 \\
36 &2\,084\,454&   3159&46&29&  & 8&  &  &  &  &  &  &  &  &  &  &   &   &   &   &   &   &    &  & 2\,087\,696 \\
37 &5\,388\,851&   9356&99&22&  & 6&  &  &  &  &  &  &  &  &  &  &   &   &   &   &   &   &    &  & 5\,398\,334 \\
38 &1\,475\,547&   6545&  &30&  & 6&  &  &  &  &  &  &  &  &  &  &   &   &   &   &   &   &    &  & 1\,482\,128 \\
39 &3\,403\,383&   7317&65&27& 2&36&  &  &  &  &  &  &  &  &  &  &   &  1&   &   &   &   &    &  & 3\,410\,831 \\
40 &   810\,399&   2084&34&48&  & 2& 4&  &  &  &  &  & 2&  &  &  &   &   &   &   &   &   &    &  &    812\,573 \\
41 &1\,645\,374&   5231&  &13&  &  &  &  &  &  &  &  &  &  &  &  &   &   &   &   &   &   &    &  & 1\,650\,618 \\
42 &   334\,536&   3308&35&39&  &28&  &  &  &  &  & 3&  &  &  & 1&  2&   &  1&   &   &   &    &  &    337\,953 \\
43 &   579\,338&   2764&32& 6&  & 6&  &  &  & 2&  &  &  &  &  &  &   &   &   &   &   &   &    &  &    582\,148 \\
44 &    94\,833&    664&  &10&  &  & 2&  &  &  &  &  &  &  &  &  &   &   &   &   &   &   &    &  &     95\,509 \\
45 &   137\,174&   1487&18&  &  &26&  &  &  &  &  & 2&  &  &  &  &   &   &   &   &   &   &    &  &    138\,707 \\
46 &    21\,818&    713& 4&  &  & 4&  &  &  &  &  &  &  &  &  &  &   &   &   &   &   &   &    &  &     22\,539 \\
47 &    18\,178&    353&  &  &  &  &  &  &  &  &  &  &  &  &  &  &   &   &   &   &   &   &    &  &     18\,531 \\
48 &       1901&    113&12& 2&  & 2& 2&  &  & 1&  &  &  &  &  &  &   &   &   &   &   &   &    &  &        2033 \\
49 &       1275&    174& 3&  & 1& 4&  &  & 2&  &  &  &  &  &  &  &   &   &   &   &   &   &    &  &        1459 \\
50 &        392&     80&  & 4&  &  &  &  &  &  &  &  &  &  &  &  &   &   &   &   &   &   &    &  &         476 \\
54 &          4&      9&  & 4&  &  &  &  &  &  &  &  &  &  &  &  &   &   &   &   &   &   &    &  &          17 \\
60 &           &       &  &  &  & 2&  &  &  &  &  &  &  &  & 2&  &   &   &   &   &  1&  1&    & 1&           7 \\
\hline
All& 55\,865\,753 & 137\,782 & 918 & 929 & 4 & 369 & 28 & 4 & 4 & 13 & 3 & 30 & 4 & 6 & 4 & 13 & 4 & 1 & 2 & 1 & 1 & 1 & 1 & 1 & 56\,005\,876 \\
\hline
\end{tabular}
\end{table*}

\[\small
\setlength{\arraycolsep}{1.5pt}
C^{n=11}_{\left|\Aut\right|=47\,520} = 
\left(
\begin{array}{ccccccccccc}
\omega&1&$-$&1&1&$-$&1&1&1&0&0\\
1&\omega&1&0&1&$-$&$-$&$-$&1&0&1\\
$-$&1&\omega&0&$-$&0&$-$&1&1&1&1\\
1&0&0&\omega&0&1&1&$-$&$-$&0&1\\
1&1&$-$&0&\omega&0&0&0&1&$-$&0\\
$-$&$-$&0&1&0&\omega&0&0&1&$-$&0\\
1&$-$&$-$&1&0&0&\omega&0&$-$&0&0\\
1&$-$&1&$-$&0&0&0&\omega&0&$-$&0\\
1&1&1&$-$&1&1&$-$&0&\omega&$-$&$-$\\
0&0&1&0&$-$&$-$&0&$-$&$-$&\omega&$-$\\
0&1&1&1&0&0&0&0&$-$&$-$&\omega
\end{array}
\right)
\]

\[\small
\setlength{\arraycolsep}{1.5pt}
C^{n=11}_{\left|\Aut\right|=1440} = 
\left(
\begin{array}{ccccccccccc}
\omega&1&$-$&1&$-$&0&0&0&0&0&1\\
1&\omega&$-$&1&1&$-$&$-$&0&1&1&0\\
$-$&$-$&\omega&1&0&1&$-$&$-$&$-$&0&$-$\\
1&1&1&\omega&0&$-$&$-$&$-$&0&$-$&1\\
$-$&1&0&0&\omega&1&$-$&1&1&$-$&1\\
0&$-$&1&$-$&1&\omega&$-$&1&$-$&1&$-$\\
0&$-$&$-$&$-$&$-$&$-$&\omega&0&0&0&0\\
0&0&$-$&$-$&1&1&0&\omega&0&0&1\\
0&1&$-$&0&1&$-$&0&0&\omega&0&$-$\\
0&1&0&$-$&$-$&1&0&0&0&\omega&$-$\\
1&0&$-$&1&1&$-$&0&1&$-$&$-$&\omega
\end{array}
\right)
\]

\[\small
\setlength{\arraycolsep}{1.5pt}
C^{n=11}_{\alpha=0} = 
\left(
\begin{array}{ccccccccccc}
\omega&1&0&$-$&1&0&$-$&1&0&0&0\\
1&\omega&1&0&$-$&$-$&1&0&0&0&0\\
0&1&\omega&1&$-$&1&0&0&0&1&0\\
$-$&0&1&\omega&$-$&$-$&0&$-$&0&$-$&$-$\\
1&$-$&$-$&$-$&\omega&$-$&$-$&0&0&0&1\\
0&$-$&1&$-$&$-$&\omega&$-$&$-$&$-$&1&0\\
$-$&1&0&0&$-$&$-$&\omega&0&$-$&0&1\\
1&0&0&$-$&0&$-$&0&\omega&$-$&$-$&$-$\\
0&0&0&0&0&$-$&$-$&$-$&\omega&1&1\\
0&0&1&$-$&0&1&0&$-$&1&\omega&0\\
0&0&0&$-$&1&0&1&$-$&1&0&\omega
\end{array}
\right)
\]

We find that there are 6493 optimal $(12,3^{12},6)$ codes with 27
distinct weight enumerators and 32 different values for
$\left|\Aut(\mathcal{C})\right|$. The weight enumerators that exist
are given by $W_{12,\alpha}(y) = 1 + (480+4\alpha)y^6 +
(3456-24\alpha)y^7 + (15120+60\alpha)y^8 + (55520-80\alpha)y^9 +
(133920+60\alpha)y^{10} + (19536-24\alpha)y^{11} +
(129408+4\alpha)y^{12}$ for all integer values $\alpha \in \{0$, $1$,
$3$, $4$, $7$, $9$, $12$, $13$, $16$, $19$, $21$, $25$, $27$, $28$,
$31$, $36$, $37$, $39$, $43$, $48$, $49$, $52$, $57$, $63$, $64$,
$81$, $144\}$.  Table~\ref{tab:opt12} gives the number of
$(12,3^{12},6)$ codes for each possible weight enumerator and
automorphism group size.  Generator matrices for all the optimal codes
of length 12 can be obtained from
\url{http://www.ii.uib.no/~larsed/nonbinary/}.  We here list generator
matrices for the unique code with maximal automorphism group size
(2\,280\,960) and a code with weight enumerator $W_{12,0}(y)$. In the
latter case, we choose the single code with maximal number of
automorphisms (11\,520).

\begin{table*}
\renewcommand{\arraystretch}{1.3}
\setlength{\tabcolsep}{3.1pt}
\setlength{\extrarowheight}{1pt}
\caption{Number of $(12,3^{12},6)$ Codes with $\left|\Aut(\mathcal{C})\right| = \beta$ and Weight Enumerator~$W_{12,\alpha}(y)$}%
\label{tab:opt12}
\centering
\begin{tabular}{crrrrrrrrrrrrrrrrrrrrrrrrrrrr}
\hline
$\beta \backslash \alpha$ & 
0 & 1 & 3 & 4 & 7 & 9 & 12 & 13 & 16 & 19 & 21 & 25 & 27 & 28 & 31 & 36 & 37 & 39 & 43 & 48 & 49 & 52 & 57 & 63 & 64 & 81 & 144 & All  \\
\hline
2           & 55 & 117 &  54 & 120 & 186 & 209 & 158 & 338 & 325 & 448 & 418 & 236 & 160 & 268 & 162 & 57 & 86 & 36 &  6 &   &  6 &   &  &  &  & & & 3445 \\
4           & 62 & 135 & 102 & 147 & 184 &  85 & 124 & 214 & 161 & 188 & 222 & 113 &  90 & 134 &  88 & 56 & 48 & 36 & 24 & 7 & 16 &   &  &  &  & & & 2236 \\
6           &    &   5 &   2 &   6 &   2 &   6 &     &     &   4 &     &   8 &   2 &   4 &     &   2 &  4 &    &    &    &   &  1 &   &  &  &  & & & 46 \\
8           & 23 &  25 &   8 &  31 &  40 &  16 &  32 &  28 &  28 &  22 &  18 &  14 &  14 &  26 &  16 & 15 &  6 &  4 &    &   &  8 & 2 &  &  &  & & & 376 \\
12          & 1 & 6 & 12 & 10 & 4 & 19 & 12 & 8 & 7 & 6 & 12 & 1 & 6 & 2 & 6 & 6 & & 2 & 2 & & 2 & 2 & & & & & & 126 \\
16          & 20 & & 2 & 4 & 6 & 2 & 5 & 4 & 18 & 4 & 4 & 2 & 2 & 8 & & 3 & 2 & 2 & & 4 & & & & & & & & 92 \\
24          & 2 & 4 & 7 & 6 & & 8 & & & 1 & & 8 & & 2 & 4 & 2 & 4 & & 8 & 2 & 2 & & & 4 & & & & & 64 \\
32          & 5 & & & 3 & & & 1 & & 5 & & & & & 2 & & & & & & & & & & & & & & 16 \\
36          & & & & & & 1 & & & & & & & 2 & & & & & & & & & & & & & & & 3 \\
48          & 3 & & & 1 & 2 & 5 & 2 & & 1 & & 6 & & & 2 & & 4 & & & & 2 & 2 & & & 2 & & & & 32 \\
64          & 2 & & & 4 & & & & & 4 & & & & & & & & & & & 1 & & 2 & & & 1 & & & 14 \\
72          & 2 & & & & & & & & & & & & 2 & & & & & & & & & & & & & & & 4 \\
80          & 1 & & & & & & & & & & & & & & & & & & & & & & & & & & & 1 \\
96          & 3 & & & & & & & & 1 & & & & & & & 1 & & & & & & & & & & & & 5 \\
108         & & & & & & & & & & & & & 2 & & & & & & & & & & & & & & & 2 \\
120         & & & & & & 1 & & & & & & & & & & & & & & & & & & & & & & 1 \\
144         & 2 & & & & & & & & & & & & 6 & & & & & & & & & & & 2 & & & & 10 \\
192         & & & & & & & & & 1 & & & & & & & 1 & & & & & & & & & & & & 2 \\
216         & 1 & & & & & 1 & & & & & & & & & & & & & & & & & & & & & & 2 \\
256         & 1 & & & & & & & & & & & & & & & & & & & & & & & & & & & 1 \\
288         & 1 & & & & & & & & & & & & & & & & & & & & & & & & & & & 1 \\
576         & 1 & & & & & & & & & & & & & & & & & & & & & & & & & & & 1 \\
720         & & & & & & & & & & & & & & & & 1 & & & & & & & & & & & & 1 \\
768         & 1 & & & & & & & & & & & & & & & & & & & & & & & & & & & 1 \\
960         & & & & & & & & & & & & & & & & & & & & & & & & & 1 & & & 1 \\
1296        & & & & & & & & & & & & & & & & & & & & & & & & & & 2 & & 2 \\
1536        & & & & & & & & & 1 & & & & & & & & & & & & & & & & & & & 1 \\
1728        & & & & & & & & & & & & & & & & 1 & & & & & & & & & & & & 1 \\
2592        & & & & & & & & & & & & & & & & & & & & & & & & 2 & & & & 2 \\
4320        & & & & & & & & & & & & & & & & 1 & & & & & & & & & & 1 & & 2 \\
11\,520     & 1 & & & & & & & & & & & & & & & & & & & & & & & & & & & 1 \\
2\,280\,960 & & & & & & & & & & & & & & & & & & & & & & & & & & & 1 & 1 \\
\hline
All&  187 & 292 & 187 & 332 & 424 & 353 & 334 & 592 & 557 & 668 & 696 & 368 & 290 & 446 & 276 & 154 & 142 & 88 & 34 & 16 & 35 & 6 & 4 & 6 & 2 & 3 & 1 & 6493 \\
\hline
\end{tabular}
\end{table*}

\[\small
\setlength{\arraycolsep}{1.5pt}
C^{n=12}_{\left|\Aut\right|=2\,280\,960} = 
\left(
\begin{array}{cccccccccccc}
\omega&$-$&$-$&$-$&$-$&$-$&1&0&1&$-$&$-$&0\\
$-$&\omega&$-$&$-$&$-$&$-$&1&$-$&$-$&0&1&0\\
$-$&$-$&\omega&$-$&$-$&$-$&$-$&$-$&1&1&0&0\\
$-$&$-$&$-$&\omega&$-$&$-$&0&1&$-$&1&$-$&0\\
$-$&$-$&$-$&$-$&\omega&$-$&$-$&1&0&$-$&1&0\\
$-$&$-$&$-$&$-$&$-$&\omega&0&0&0&0&0&0\\
1&1&$-$&0&$-$&0&\omega&1&1&1&1&$-$\\
0&$-$&$-$&1&1&0&1&\omega&1&1&1&$-$\\
1&$-$&1&$-$&0&0&1&1&\omega&1&1&$-$\\
$-$&0&1&1&$-$&0&1&1&1&\omega&1&$-$\\
$-$&1&0&$-$&1&0&1&1&1&1&\omega&$-$\\
0&0&0&0&0&0&$-$&$-$&$-$&$-$&$-$&\omega
\end{array}
\right)
\]

\[\small
\setlength{\arraycolsep}{1.5pt}
C^{n=12}_{\alpha=0,\, \left|\Aut\right|=11\,520} = 
\left(
\begin{array}{cccccccccccc}
\omega&$-$&$-$&$-$&$-$&$-$&1&$-$&$-$&1&0&0\\
$-$&\omega&$-$&$-$&$-$&$-$&0&$-$&0&1&$-$&1\\
$-$&$-$&\omega&$-$&$-$&$-$&0&0&$-$&1&1&$-$\\
$-$&$-$&$-$&\omega&$-$&$-$&1&0&$-$&0&$-$&1\\
$-$&$-$&$-$&$-$&\omega&$-$&1&$-$&0&0&1&$-$\\
$-$&$-$&$-$&$-$&$-$&\omega&0&0&0&0&0&0\\
1&0&0&1&1&0&\omega&$-$&$-$&$-$&1&1\\
$-$&$-$&0&0&$-$&0&$-$&\omega&$-$&$-$&1&1\\
$-$&0&$-$&$-$&0&0&$-$&$-$&\omega&$-$&1&1\\
1&1&1&0&0&0&$-$&$-$&$-$&\omega&1&1\\
0&$-$&1&$-$&1&0&1&1&1&1&\omega&1\\
0&1&$-$&1&$-$&0&1&1&1&1&1&\omega
\end{array}
\right)
\]

\section{Conclusion}\label{sec:conc}

According to the mass formula bound (\ref{eq:massbound}), the total
number of codes of length 11 and 12 are $t_{11} \ge 1\,592\,385\,579$
and $t_{12}\ge 2\,938\,404\,780\,748$, which makes complete
classifications infeasible, at least with our computational
resources. Running our algorithm on a typical desktop computer, the
classification of codes of length $n$ was completed in less than five
minutes for $n \le 8$, about two hours for $n=9$, and about a week for
$n=10$. Most of this time is spent canonizing the equivalence graphs
with \emph{nauty}, and far more time is used on codes with large
automorphism groups than on codes with trivial or small automorphism
groups.  This means that our previous classification
algorithm~\cite{nonbinary}, using local complementation, might still
be useful in some cases, since we observe that graphs corresponding to
codes with large automorphism groups typically have small LC orbits.
For instance, we could speed up our classification algorithm by not
only removing isomorphisms from the set of lengthened codes, but also
generating and storing a limited number of LC orbit members of each
graph, and checking new graphs for isomorphism against this set.
Finding all optimal codes of length 11 and 12 required 80 and 320 days
of CPU time, respectively, and a parallel cluster computer was used
for this search.  We observed that most of this time was spent on
computing minimum distance to eliminate non-optimal codes, and much
less time on canonizing the optimal codes.

Although this paper has focused on codes over $\GF(9)$, our
classification algorithm can be generalized to Hermitian self-dual
additive codes over $\GF(q=m^2)$ for any prime power $m$. (One simply
needs to find an appropriate coordinate graph, as discussed in
Section~\ref{sec:eqgraph}.) The results in this paper also has
applications beyond the study of additive codes. The correspondence
between self-dual additive codes over $\GF(9)$ and 3-weighted graphs
means that we have also classified particular classes of 3-weighted
graphs that should be of interest in graph theory. An equivalence
class of self-dual additive codes over $\GF(9)$ maps to an orbit of
graphs under generalized local complementation~\cite{newlc,nonbinary}.
Orbits of graphs with respect to local complementation has a long
history in combinatorics~\cite{hubert, bouchet}, with several
applications, for instance in the theory of \emph{interlace
  polynomials}~\cite{aigner, interlacecodes}. The generalization to
weighted graphs is a natural next step.  The results in this paper
also have applications in the field of quantum information theory.
Our previous classification of codes over $\GF(4)$~\cite{selfdualgf4}
has since led to new results in the study of the \emph{entanglement}
of \emph{quantum graph states}~\cite{cabello}, and the new data
obtained in this paper will yield similar insights into the properties
of ternary quantum graph states.

\section*{Acknowledgement}
Thanks to Parallab, the High Performance Computing unit of the
University of Bergen, whose cluster computer made the computational
results in this paper possible.
The author would like to thank the anonymous reviewers for providing
useful suggestions and corrections that improved the quality of the
manuscript.

\IEEEtriggeratref{14}

\end{document}